\documentclass{amsproc}
\usepackage{amsmath,amsthm,amsfonts,rotating}

\newtheorem{theorem}{Theorem}[section]
\theoremstyle{definition}
\newtheorem{definition}[theorem]{Definition}
\newtheorem{proposition}[theorem]{Proposition}
\newtheorem{corollary}[theorem]{Corollary}
\newtheorem{example}[theorem]{Example}

\numberwithin{equation}{section}

\begin{document}

\title[Factorization theory for Wiener-Hopf plus Hankel operators]{Factorization theory
for Wiener-Hopf plus Hankel operators with almost periodic symbols}

\author{A.~P.~Nolasco}
\address{Department of Mathematics, University of Aveiro,
3810-193 Aveiro, Portugal}
\email{anolasco@mat.ua.pt}
\thanks{A.~P.~Nolasco is sponsored by {\em Funda\c{c}\~{a}o para a Ci\^{e}ncia e a Tecnologia}
(Portugal) under grant number SFRH/BD/11090/2002.}

\author{L.$\!$ P.~Castro}
\address{Department of Mathematics, University of Aveiro,
3810-193 Aveiro, Portugal} \email{lcastro@mat.ua.pt}
\urladdr{http://www.mat.ua.pt/lcastro/}




\maketitle

\begin{abstract}
A factorization theory is proposed for Wiener-Hopf plus Hankel
operators with almost periodic Fourier symbols. We introduce a
factorization concept for the almost periodic Fourier symbols such
that the properties of the factors will allow corresponding
operator factorizations. Conditions for left, right, or both-sided
invertibility of the Wiener-Hopf plus Hankel operators are
therefore obtained upon certain indices of the factorizations.
Under such conditions, the one-sided and two-sided inverses of the
operators in study are also obtained.
\end{abstract}

\section{Introduction}
\subsection{Motivation and Historical Background}

Wiener-Hopf plus Hankel operators (as well as their discrete
analogues based on Toeplitz and Hankel operators) are well known
to play an important role in several applied areas. E.g., this is
the case in certain wave diffraction problems, digital signal
processing, discrete inverse scattering, and linear prediction.
For concrete examples of a detailed enrollment of those operators
in these (and others) applications we refer to \cite{CaSpTe03,
CaSpTe04otaa, FaYa92, LeMeTe92, MeTh92, MeSpTe92, SaLe93, Te91}.

In view of the needs of the applications, it is natural to expect
a corresponding higher interest in mathematical fundamental
research for those kind of operators. In fact, in recent years
several authors have contributed to the mathematical understanding
of Wiener-Hopf plus Hankel (and their discrete analogues) under
different types of assumptions (cf.~\cite{BaTo04, BaEh05,
BaEhWi03, CaSp05, CaSpTe04mn, Eh04, He02, HeRo05, MeSpTe92,
NoCa04, RoSaSi97, RoSi90, Si87}).

We recall that the name ``Wiener-Hopf operators" is due to the
initial work of Norbert Wiener and Eberhard Hopf \cite{WeHo31} where
 a reasoning to solve integral equations whose
kernels depend only on the difference of the arguments was
provided:
\begin{equation}\label{WH:Eq}
c f(x) + \int_0^{+\infty} k(x-y)f(y)dy=g(x), \quad x \in
\mathbb{R}_+,
\end{equation}
i.e.~the so-called {\em integral Wiener-Hopf equations}. Here $c \in
\mathbb{C},$ $k \in L^1(\mathbb{R})$ and $f,g \in
L^2(\mathbb{R}_+)$, where $c$ and $k$ are fixed, $g$ is given and
$f$ is the unknown element.

From those Wiener-Hopf equations arise the {\em (classical)
Wiener-Hopf operators} defined by
\begin{equation}\label{WH:Op}
W_\phi f(x)= c f(x) + \int_0^{+\infty} k(x-y)f(y)dy \;\;,\quad x \in
\mathbb{R}_+,
\end{equation}
where $\phi$ belongs to the {\em Wiener algebra}:
\begin{equation*}
\mathbb{W}=\{\phi:\phi=c+\mathcal{F}k,c \in \mathbb{C},k \in
L^1(\mathbb{R})\}
\end{equation*}
(which is a Banach algebra when endowed with the norm
$\|c+\mathcal{F}k\|_\mathbb{W}=|c|+\|k\|_{L^1(\mathbb{R})}$ and
the usual multiplication operation). Having in mind the
convolution operation, the definition of $W_\phi$ in (\ref{WH:Op})
gives rise to an understanding of the Wiener-Hopf operators as
convolution type operators. Therefore, they can also be
represented as
\begin{equation}\label{WH:geral}
W_\phi = r_+ \mathcal{F}^{-1}\phi\cdot\mathcal{F} \;:\;
L^2_+(\mathbb{R})\rightarrow
L^2(\mathbb{R_+}).
\end{equation}
Here $L^2_+(\mathbb{R})$ denotes the subspace of $L^2(\mathbb{R})$
formed by all the functions supported in the closure of
$\mathbb{R}_+=(0,+\infty)$, $r_+$ is the {\em restriction operator}
from $L^2(\mathbb{R})$ into $L^2(\mathbb{R_+})$, and $\mathcal{F}$
denotes the {\em Fourier transformation}. The Wiener-Hopf operators
on $L^2_+(\mathbb{R})$ may also be written in the form
\begin{eqnarray} \label{W1}
P_+ A = \ell_0 r_+
\mathcal{F}^{-1}\phi\cdot\mathcal{F}:L^2_+(\mathbb{R})\rightarrow
L^2_+(\mathbb{R})\,,
\end{eqnarray}
where $\ell_0 : L^2(\mathbb{R_+})\rightarrow L^2_+(\mathbb{R})$ is the {\em
zero extension operator}, $P_+=\ell_0 r_+$ is the {\em canonical projection}
of $L^2(\mathbb{R})$ onto $L^2_+(\mathbb{R})$, and $A$ is the translation
invariant operator $\mathcal{F}^{-1}\phi\cdot\mathcal{F}$. Looking now to
the structure of the operators in (\ref{WH:geral}) and (\ref{W1}), we
recognize that possibilities other than only the Wiener algebra can be
considered for the so-called {\em Fourier symbols} $\phi$ of the Wiener-Hopf
operators. Namely, we may consider to choose $\phi$ among the
$L^{\infty}(\mathbb{R})$ elements (i.e., as a measurable essentially bounded
function in the real line).

Within the context of (\ref{WH:Eq}) and (\ref{WH:Op}), the {\em
Hankel integral operators} have the form
\begin{equation}\label{H:Op}
H f(x)= \int_0^{+\infty} k(x+y)f(y)dy \;\;,\quad x \in \mathbb{R}_+
\end{equation}
(for some $k \in L^1(\mathbb{R})$). It is well known that $H$, as
an operator defined between $L^2$ spaces, is a compact operator.
However, as seen above, it is also possible to provide a rigorous
meaning to the expression (\ref{H:Op}) when the kernel $k$ is a
temperate distribution whose Fourier transform belongs to
$L^{\infty}(\mathbb{R})$.

We would like to mention that the discrete analogue of $H$ has its
roots in the year of 1861 with the Ph.D.~thesis of Hermann
Hankel~\cite{Ha61}. There the study of finite matrices with
entries depending only on the sum of the coordinates was proposed.
Determinants of infinite complex matrices with entries defined by
$a_{jk}=a_{j+k}$ (for $j,k \geq 0$, and where $a=\{a_j\}_{j \geq
0}$ is a sequence of complex numbers) were also studied. For these
(infinite) Hankel matrices, one of the first main results was
obtained by Kronecker in 1881 \cite{Kr81} when characterizing the
Hankel matrices of finite rank as the ones that have corresponding
power series, $a(z)=\sum_{j=0}^\infty a_j\, z^j$, which are
rational functions.
 In 1906, Hilbert proved that the operator (induced by the famous Hilbert
matrix), $\mathcal{H}:\ell^2\rightarrow\ell^2\;, \;\{b_j\}_{j \geq
0} \mapsto \left\{\sum_{k=0}^\infty {b_k}/{(j+k+1)}  \right\}_{j
\geq 0} \,$, is bounded on $\ell^2$. This result may be viewed as
the origin of (discrete) Hankel operators, as natural objects
arising from Hankel matrices. Later on, in 1957, Nehari presented a
characterization of bounded Hankel operators on $\ell^2$
\cite{Ne57}. Due to the importance of such characterization, we may
say that it marks the beginning of the contemporary period of the
study of Hankel operators.

As for combinations between Wiener-Hopf (or Toeplitz) and Hankel
operators, an important initial step occurred in 1979 when Power
\cite{Po79} used the $C^\ast$--algebra generated by Toeplitz and
Hankel operators. In particular, Power devoted a particular
attention to those kinds of operators having piecewise continuous
Fourier symbols. Later on, several other authors considered also
interactions between Wiener-Hopf and Hankel operators, as well as
the algebra generated by them (see e.g.~\cite{BaTo04, BaEh05,
BaEhWi03, BoSi90, Po80, Po82, RoSi90, Si87}). As a consequence,
the theory of Wiener-Hopf plus Hankel operators is nowadays well
developed for some classes of Fourier symbols (like in the case of
continuous or piecewise continuous symbols). In particular, the
invertibility and Fredholm properties of such kind of operators
with piecewise continuous Fourier symbols are now well known (and
of great importance for the applications \cite{CaSpTe04otaa,
LeMeTe92, MeSpTe92, Te91}). However, this is not the case for
almost periodic Fourier symbols which are also important in the
applications in view of their appearance due to, e.g.,
(i)~particular finite boundaries in the geometry of physical
problems \cite{CaSpTe03}, or (ii)~the needs of compositions with
shift operators which introduce almost periodic elements in the
Fourier symbols of those operators \cite{CaSpIEOT00, KaSp94}.

The present paper introduces characterizations of the
invertibility and Fredholm properties of Wiener-Hopf plus Hankel
operators with almost periodic Fourier symbols.

\subsection{The Wiener-Hopf plus Hankel Operators in Study, Basic Definitions, and Main Results Outline}
The main objects of the present work are the
{Wiener-Hopf plus Hankel operators} with Fourier symbols in the
algebra of almost periodic functions, and acting between $L^2$
Lebesgue spaces. In a detailed way, we will consider operators with
the form
\begin{equation}\label{opWH1}
W\!H_\phi \,=\, W_\phi+H_\phi \;:\; L^2_+(\mathbb{R})\rightarrow
L^2(\mathbb{R_+}),
\end{equation}
with $W_\phi$ and $H_\phi$ being {\em Wiener-Hopf} and {\em Hankel
operators} defined by
\begin{eqnarray}\label{op WH}
W_\phi &\!\!\!=\!\!\!& r_+ \mathcal{F}^{-1}\phi\cdot\mathcal{F}
\;:\; L^2_+(\mathbb{R})\rightarrow
L^2(\mathbb{R_+})\\
\label{op H} H_\phi &\!\!\!=\!\!\!& r_+
\mathcal{F}^{-1}\phi\cdot\mathcal{F}J \;:\;
L^2_+(\mathbb{R})\rightarrow L^2(\mathbb{R_+}),
\end{eqnarray}
respectively. Here and in what follows, $J$ is the {\em reflection
operator} given by the rule $J\varphi(x)=\widetilde{\varphi}(x)=
\varphi(-x), \; x \in \mathbb{R}$. According to \eqref{opWH1},
\eqref{op WH} and \eqref{op H}, we have
\begin{equation*}
W\!H_\phi= r_+ (\mathcal{F}^{-1}\phi\cdot\mathcal{F} +
\mathcal{F}^{-1}\phi\cdot\mathcal{F}J)\\
          = r_+ \mathcal{F}^{-1}\phi\cdot\mathcal{F} (I_{L^2_+(\mathbb{R})} +
          J),
\end{equation*}
where $I_{L^2_+(\mathbb{R})}$ denotes the {\em identity operator} in
$L^2_+(\mathbb{R})$. Furthermore, since
$$I_{L^2_+(\mathbb{R})} +
J =\ell^e r_+\,,$$ where $\ell^e : L^2(\mathbb{R_+}) \rightarrow
L^2(\mathbb{R})$ denotes the {\em even extension operator}, we may write the
Wiener-Hopf plus Hankel operator as
\begin{equation}\label{eq17}
W\!H_\phi= r_+ \mathcal{F}^{-1}\phi\cdot\mathcal{F} \ell^e r_+\,.
\end{equation}
Finally, the Fourier symbol $\phi$ belongs to the algebra $AP$ of
{\em almost periodic functions}, i.e., the smallest closed
subalgebra of $L^\infty(\mathbb{R})$ that contains all the functions
$e_\lambda \; (\lambda \in \mathbb{R})$ with
$e_{\lambda}(x)=e^{i\lambda x}$, $x\in \mathbb{R}$.

We will proceed now with some basic definitions which will be needed
in what follows. Let $T:X\rightarrow Y$ be a bounded linear operator
acting between Banach spaces. The operator $T$ is said to be {\em
normally solvable} if $\text{Im}\,T$ is closed. In this case, the
{\em cokernel of T} is defined as $\text{Coker}\,T=Y/ \text{Im}\,T$.
For a normally solvable operator $T$, the {\em deficiency numbers}
of $T$ are given by
\begin{equation}\label{KerCoker}
n(T):=\dim \text{Ker}\,T, \quad d(T):=\dim \text{Coker}\,T.
\end{equation}
If at least one of the deficiency numbers is finite, the operator $T$ is
said to be a \emph{semi-Fredholm} operator. A normally solvable operator $T$
is said to be: (i) a {\em Fredholm} operator if both $n(T)$ and $d(T)$ are
finite; (ii) {\em left-Fredholm} if $n(T)$ is finite; (iii) {\em right
Fredholm} if $d(T)$ is finite. In the case when only one of the deficiency
numbers is finite, the operator $T$ is said to be a \emph{properly
semi-Fredholm} operator. In a more detailed way, a normally solvable
operator $T$ is said to be {\em properly left-Fredholm} if $n(T)$ is finite
and $d(T)$ is infinite, and {\em properly right-Fredholm} if $d(T)$ is
finite and $n(T)$ is infinite. We point out that in German and Russian
literature, \mbox{(semi-)}Fredholm operators are often called
(semi-)\emph{Noether operators}. This is due to the pioneering work
\cite{No21} of Fritz Noether who was the first to discover that
 singular integral operators with nonvanishing continuous symbols
 are normally solvable, and have finite kernel and cokernel
dimensions. Once again, in German and Russian literature (and related with
the notation used in (\ref{KerCoker})), right-Fredholm operators and left
Fredholm operators are frequently called \emph{$n$-normal operators} and
\emph{$d$-normal operators}, respectively (see, e.g.,
\cite{BoKaSp02,KaSp86,KaSp90,KaSp94}).

 Let us choose the notation $\mathcal{G}B$ for the
group of all invertible elements of a Banach algebra $B$. By Bohr's
theorem, for each $\phi \in \mathcal{G}AP$ there exists a real
number $\kappa(\phi)$ and a function $\psi \in AP$ such that
\begin{eqnarray}\label{meanM}
\phi(x)=e^{i\kappa(\phi)x}e^{\psi(x)}\,, \quad x\in \mathbb{R}.
\end{eqnarray}
Since $\kappa(\phi)$ is uniquely determined, $\kappa(\phi)$ is
usually called the {\em mean motion} of $\phi$.

For Wiener-Hopf operators with Fourier symbols in $\mathcal{G}AP$,
there is a famous semi-Fredholm and invertibility criterion due to
Gohberg-Feldman/Coburn-Douglas based on the sign of the mean motion
of the Fourier symbol of the operator (cf. \cite{CoDo} and
\cite{GoFe} or \cite[Theorem 2.28]{BoKaSp02}):
\begin{itemize}
\item[(a)] if the mean motion of the Fourier symbol is negative, then the
Wiener-Hopf operator is properly right-Fredholm and right-invertible;
\item[(b)] if the mean motion of the symbol is positive, then the
Wiener-Hopf operator is properly left-Fredholm and left-invertible;
\item[(c)] if the mean motion of the symbol is zero, then
the Wiener-Hopf operator is invertible.
\end{itemize}
This criterion was one of the initial motivations for the present work.
Accordingly, the main purpose of this paper is to establish an invertibility
and Fredholm criterion for Wiener-Hopf plus Hankel operators with almost
periodic Fourier symbols. To reach such criterion, we need to introduce a
new factorization concept for $AP$ functions -- the so-called $AP$
asymmetric factorization. As we will see in Section~\ref{sect4}
(Definition~\ref{def4.1}), a function $\phi \in \mathcal{G}AP$ is said to
admit an $AP$ asymmetric factorization if it can be represented in the form
$\phi=\phi_{-} \; e_{\lambda} \; \phi_{e}$, where $\lambda \in{\mathbb{R}}$,
$\phi_{-} \in \mathcal{G}AP^{-}$ (cf.~(\ref{AP-})), and $\phi_{e} \in
\mathcal{G}L^\infty(\mathbb{R})$ with $\widetilde{\phi_e}=\phi_e$.
Therefore, assuming that $\phi \in \mathcal{G}AP$ admits an $AP$ asymmetric
factorization, the obtained invertibility and Fredholm criterion for
Wiener-Hopf plus Hankel operators $W\!H_\phi$ (cf.~Theorem~\ref{th1}) has a
similar structure as the theorem of Gohberg-Feldman/Coburn-Douglas mentioned
before, and states the following:
\begin{itemize}
\item[(a)] if $\lambda < 0$, then the Wiener-Hopf plus Hankel operator is properly
right-Fredholm and right-invertible;
\item[(b)] if $\lambda > 0$,
then the Wiener-Hopf plus Hankel operator is properly left-Fredholm and
left-invertible;
\item[(c)] if $\lambda = 0$, then
the Wiener-Hopf plus Hankel operator is invertible.
\end{itemize}
Under such conditions, we can do even better: we can provide a
formula for the one-sided and two-sided inverses of the
Wiener-Hopf plus Hankel operators by using the factors of the $AP$
asymmetric factorization. This result is stated in the last
section, in Theorem~\ref{th3}, and exhibits the importance of
having convenient factorizations for the Fourier symbols of the
corresponding operators.

We will also present a result on the invertibility dependencies between
Wiener-Hopf and Wiener-Hopf plus Hankel operators with the same $AP$ Fourier
symbol (cf.~Corollary~\ref{Cor5.3}). At a first glance, this may appear to
be a very surprising result. Noticing however that (in this case) we will be
dealing with a particular kind of factorization, it is more natural to hope
to achieve the invertibility of Wiener-Hopf plus Hankel operators from the
invertibility of Wiener-Hopf operators (by using certain relations between
different $AP$ factorizations).

Before arriving at the main Section~\ref{sect5} (where the last
briefly described results will appear in detail), we will recall in
Section~\ref{sect2} several useful details about almost periodic
functions, operator identities for Wiener-Hopf plus Hankel operators
are exhibited in Section~\ref{sect3} (allowing therefore the
understanding of certain compositions of those kind of operators),
and factorization concepts are proposed and analyzed in
Section~\ref{sect4}.

\section{Almost Periodic Functions}\label{sect2}

The theory of almost periodic functions (\emph{$AP$ functions})
was created by Harald Bohr between 1923 and 1925. Since then, many
contributions to the development of the theory of $AP$ functions
were made (namely by V.~V.~Stepanov, H.~Wyel, A.~S.~Besicovitch,
S.~Bochner, J. Von Neumann, C.~Corduneanu, and others). The
importance of $AP$ functions in problems of differential
equations, stability theory and dynamical systems potentiated the
development of the theory of these functions.

For defining $AP$ functions, Bohr used the concepts of
\emph{relative density} and \emph{translation number}
(cf.~\cite{Bo33}). A set $E \subset \mathbb{R}$ is said to be
\emph{relatively dense} if there exists a number $l>0$ such that any
interval of length $l$ contains at least one number of $E$. In
addition (considering $\phi$ being a real or complex function
defined on the real line), a number $\tau$ is called a
\emph{translation number} of $\phi$, corresponding to
$\varepsilon>0$, if
\begin{equation*}
|\phi(x+\tau)-\phi(x)|\leq \varepsilon
\end{equation*}
for all $x \in \mathbb{R}$. These two last concepts provide the
conditions to present the alternative definition of $AP$ function: a
continuous function $\phi$ defined on the real line is called
\emph{almost periodic} if for every $\varepsilon>0$ there exists a
relatively dense set of translation numbers of $\phi$ corresponding
to $\varepsilon$. That is to say, $\phi$ is called an $AP$ function
if for every $\varepsilon>0$ there exists a number $l>0$ such that
any interval of length $l$ contains at least one number $\tau$ for
which
\begin{equation*}
|\phi(x+\tau)-\phi(x)|\leq \varepsilon\,,
\end{equation*}
for all $x \in \mathbb{R}$.

Like the periodic functions, $AP$ functions can also be represented
by Fourier series. To obtain such representation, we have to
consider the mean value of an $AP$ function. By the {\em Mean Value
Theorem}, it follows that for every $\phi \in AP$ there exists the
\emph{Bohr mean value} (or the \emph{mean value}) of $\phi$:
\begin{equation*}
M(\phi) := \lim_{T \rightarrow \infty}
\frac{1}{T}\int_0^T\phi(x)\,dx.
\end{equation*}
For $\phi\in AP$ and $\lambda\in \mathbb{R}$, the function $\phi\,
e_{-\lambda}$, being the product of two AP functions, is also an
$AP$ function (recall that $e_{\lambda}(x)=e^{i\lambda x}$, $x \in
\mathbb{R}$). Therefore, there exits the mean value of $\phi\,
e_{-\lambda}$. The set
\begin{equation*}
\Omega(\phi):=\{\lambda \in \mathbb{R}:M(\phi\, e_{-\lambda}) \neq
0\}
\end{equation*}
is at most countable, and is called the \emph{Bohr-Fourier spectrum}
of $\phi$. The Fourier series associated with the function $\phi$ is
given by
\begin{equation*}
\sum_{j} \phi_j \, e_{\lambda_j}\,,
\end{equation*}
and we may write
\begin{equation*}
\phi(x)\sim\sum_{j} \phi_j \, e^{i\lambda_jx}\,,
\end{equation*}
where $\phi_j=M(\phi\, e_{-\lambda_j})$. The elements $\lambda_j$ of
the Bohr-Fourier spectrum are called the \emph{Fourier exponents},
and the corresponding mean values $\phi_j$ are called the
\emph{Fourier coefficients}.

In \cite{Bo33}, we find the \emph{Fundamental Theorem} of the theory
of $AP$ functions. The Fundamental Theorem states that the class of
almost periodic functions is identical with the closure of the class
of all \emph{trigonometric polynomials}:
\begin{equation}\label{trigpol}
p(x)=\sum_{n=1}^N \phi_n \, e^{i\lambda_nx},
\end{equation}
where $\phi_n \in \mathbb{C}$ and $\lambda_n \in \mathbb{R}$.  To
reach to the Fundamental Theorem, Bohr proved that the limit of a
uniformly convergent sequence of $AP$ functions is also an almost
periodic function. Thus, being trigonometric polynomials almost
periodic functions, it follows that the limit of a uniformly
convergent sequence of trigonometric polynomials is an $AP$
function. This means that every function belonging to the closure of
the class of all finite sums (\ref{trigpol}) is in the class of
almost periodic functions. The converse also holds. That is, every
$AP$ function is the limit of a uniformly convergent sequence of
trigonometric polynomials. This result is called the
\emph{Approximation Theorem}. More precisely, the approximation
theorem asserts that given $\phi \in AP$, for each $\varepsilon>0$,
there exists a trigonometric polynomial $p$, whose exponents are the
Fourier exponents of $\phi$, and such that
\begin{equation*}
|\phi(x)- p(x)|\leq \varepsilon
\end{equation*}
for all $x \in \mathbb{R}$.

\begin{figure}[h]
\hfil
\includegraphics[width=5.5cm]{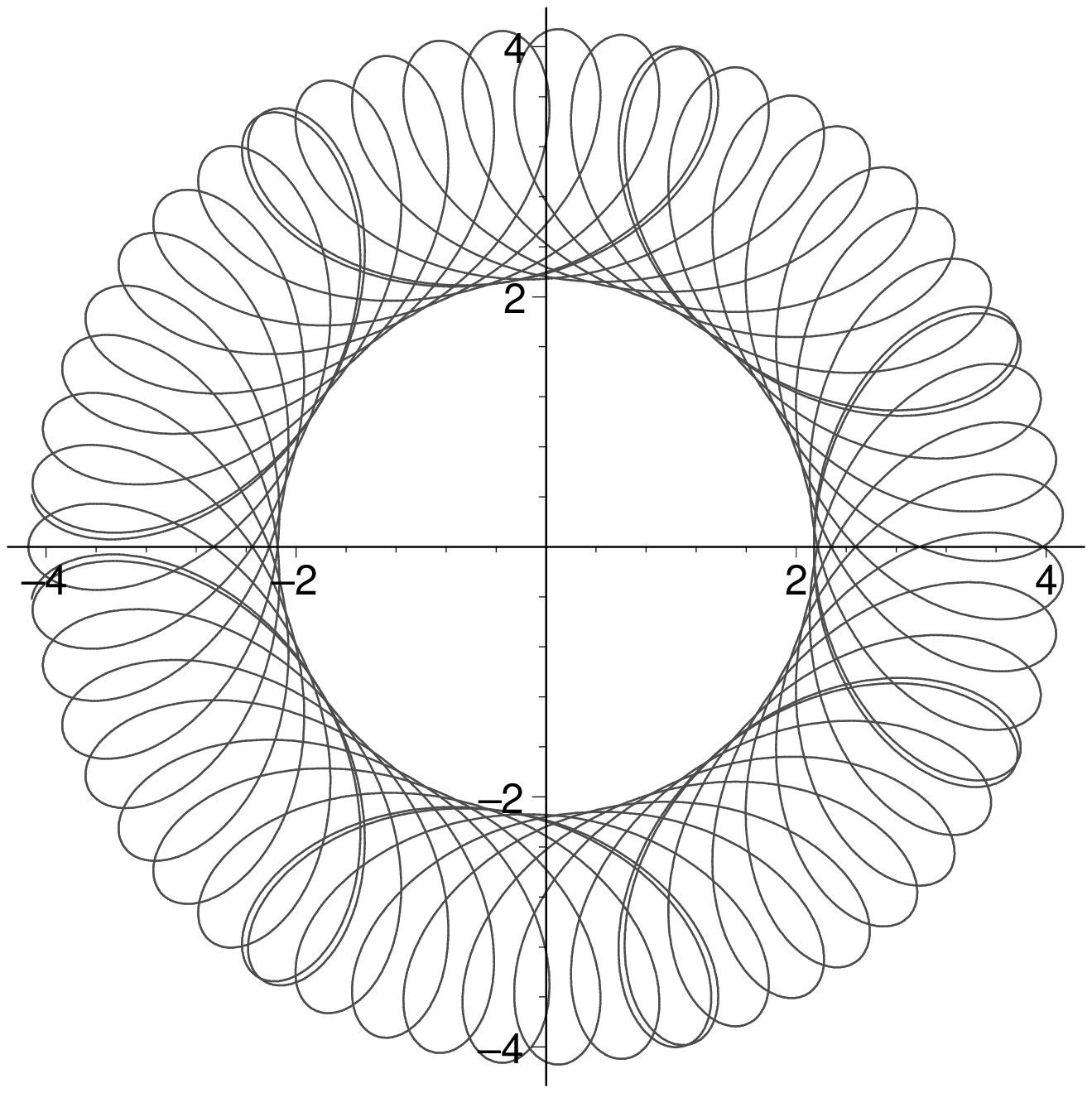} \hskip1cm\includegraphics[width=5.5cm]{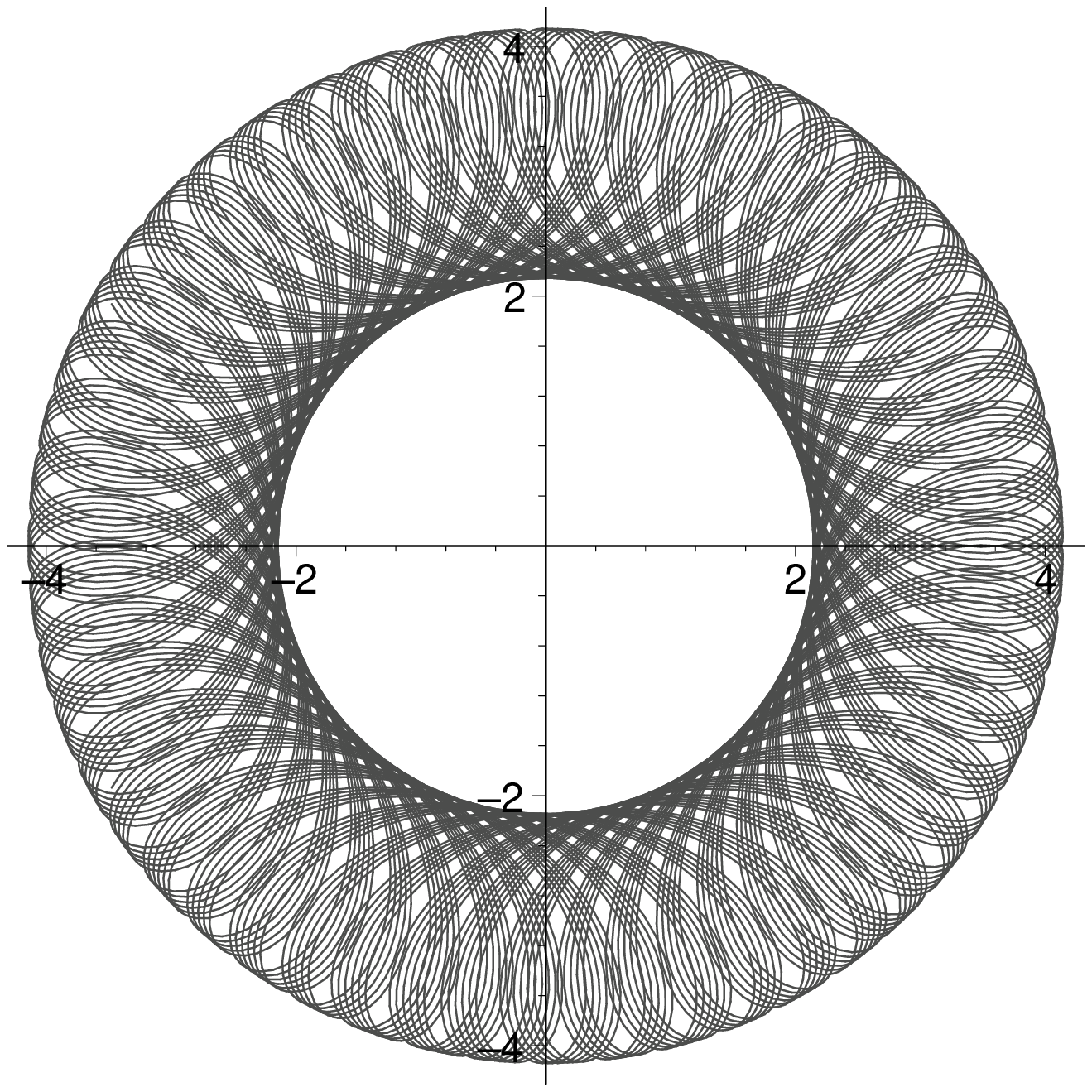}
\caption{The image of $\phi_1(x) = -e^{i \pi x}- \pi e^{-\frac{i
x}{2}}$ for $x$ in $[-50,50]$ and $[-200,200]$,
respectively.}\label{fig1}

\vspace{10mm}

\hfil
\includegraphics[width=5.5cm]{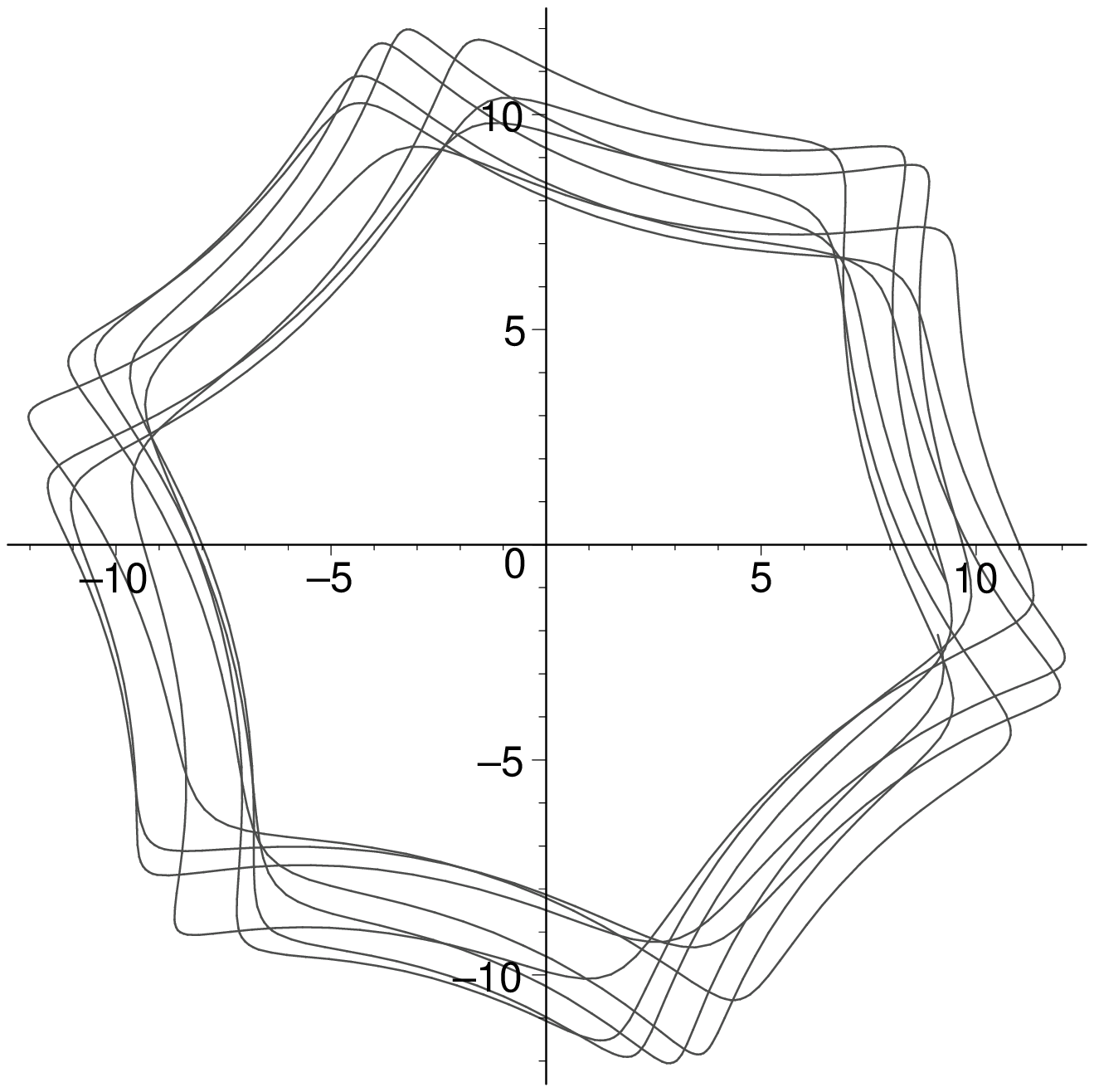}\hskip1cm\includegraphics[width=5.5cm]{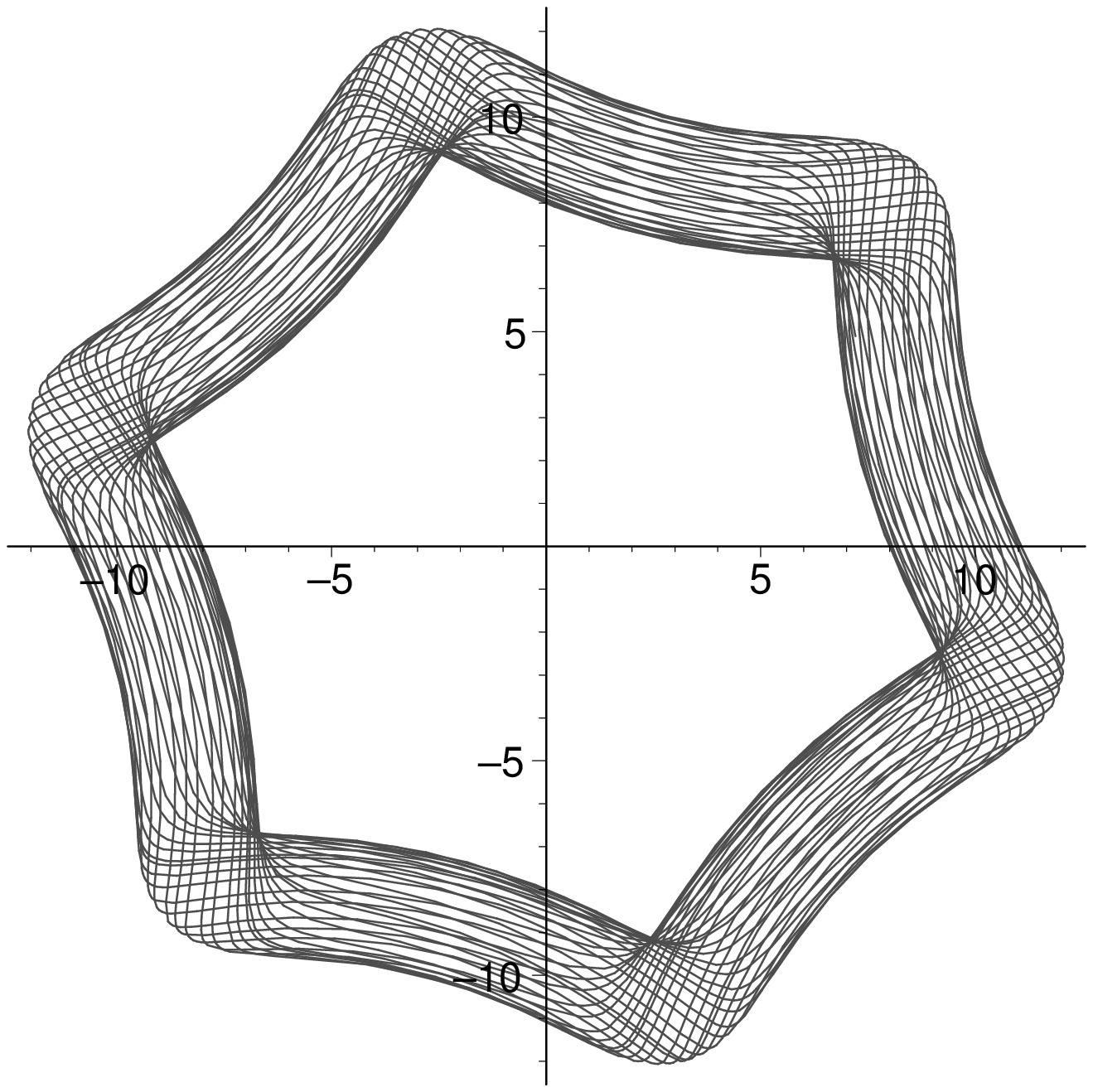}
\caption{The image of $\phi_2(x)= 10e^{-ix}+(1+i)e^{i \sqrt{2}
x}-i e^{i5x}$ for $x$ in $[-25,25]$ and $[-100,100]$,
respectively.}\label{fig2}
\end{figure}

From the Approximation Theorem, we obtain a characterization of almost
periodic functions which is sometimes presented as the definition of almost
periodic functions. See, e.g., \cite{Co} where $AP$ functions are those
complex-valued functions defined on the real line which can be uniformly
approximated by trigonometric polynomials, and \cite{BoKaSp02} where the
algebra of the $AP$ functions is defined as the smallest closed subalgebra
of $L^\infty(\mathbb{R})$ that contains all the functions $e_\lambda \;
(\lambda \in \mathbb{R})$.

To illustrate the notion of $AP$ function, Figures~\ref{fig1} and
\ref{fig2} exhibit the images of two particular $AP$ functions
$\phi_1$ and $\phi_2$.

From Bohr's Theorem (cf.~also (\ref{meanM})), it is possible to
observe that the argument of an invertible $AP$ function is given
by the sum of a linear function and an $AP$ function. For $\phi
\in \mathcal{G}AP$, the mean motion of $\phi$ can be obtained by
\begin{equation*}
\kappa(\phi)=\lim_{T\rightarrow \infty}\frac{(\arg \phi)(T)-(\arg
\phi)(-T)}{2T}\,,
\end{equation*}
where $\arg\phi$ is any continuous argument of $\phi$.

Considering that every function in $AP$ may be represented by a series, but
not every function in $AP$ may be represented by an absolutely convergent
series, it is also useful to consider the subclass $APW$  of all functions
$\varphi \in AP$ which can be written in the form of an absolutely
convergent series:
\begin{equation*}
\varphi(x)=\sum_j \varphi_j e^{i\lambda_jx}\quad (x\in
\mathbb{R}),\qquad \lambda_j \in \mathbb{R},\qquad \sum_j
|\varphi_j|<\infty\,.
\end{equation*}

To end this section about $AP$ functions, we present the definitions of
other special subsets of $AP$ that will have a preponderant use in the
proposed factorizations of $AP$ functions in below. Let $AP^-$ ($AP^+$)
denote the smallest closed subalgebra of $L^\infty(\mathbb{R})$ that
contains all the functions $e_\lambda$ with $\lambda \leq 0$ ($\lambda \geq
0$). Equivalently, we have
\begin{equation}\label{AP-}
AP^-=\{\phi \in AP:\Omega(\phi)\subset(\infty,0]\}
\end{equation}
and
\begin{equation*}
AP^+=\{\phi \in AP:\Omega(\phi)\subset [0,\infty)\}.
\end{equation*}
For more information about $AP$ functions we refer to
\cite{Be32,Bo33,Co,LeZh82}.

\section{Operator Identities for Wiener-Hopf plus Hankel Operators}\label{sect3}

The study of relations between different classes of operators is
an important subject in Operator Theory since (in particular) it
allows the transfer of properties from one class of operators to
other ones. An example of these kind of relations is the {\em
equivalence} between operators: considering two bounded linear
operators acting between Banach spaces, $T:X_1 \rightarrow X_2$
and $S:Y_1 \rightarrow Y_2$, the operators $T$ and $S$ are said to
be {\em equivalent} \cite{BaTs92} if there are two boundedly
invertible linear operators, $E:Y_2 \rightarrow X_2$ and $F:X_1
\rightarrow Y_1$, such that
\begin{equation}\label{eq10}
T=E\;S\;F.
\end{equation}
Another important example of relation between operators is the {\em
equivalence after extension} relation (or {\em matricial coupling}), which
is a generalization of the equivalence relation between operators. In
\cite{BaTs92}, Bart and Tsekanovskii proved that two operators acting
between Banach spaces are equivalent after extension if and only if they are
matricially coupled. It follows from \eqref{eq10} that if two operators are
equivalent, then they belong to the same {\em regularity class} \cite{Ca03,
CaSp98zaa, Sp85}. In particular, this implies that one of these operators is
invertible, one-sided invertible, Fredholm, properly left-Fredholm, properly
right Fredholm or normally solvable, if and only if the other operator
enjoys the same property. In our case we are interested in obtaining
relations between Wiener-Hopf plus Hankel operators and Wiener-Hopf
operators, and this idea will be used ahead in the process of obtaining an
invertibility criterion for Wiener-Hopf plus Hankel operators.

Therefore, in order to reach such relation, we will start by
identifying certain operator identities for Wiener-Hopf plus Hankel
operators and Wiener-Hopf operators.

From the Wiener-Hopf and Hankel operator theory, the following
relations are well known:
\begin{align*}
&W_{\phi\varphi}=W_\phi\,\ell_0\,W_\varphi+H_\phi\,\ell_0\,H_{\widetilde{\varphi}},\\
&H_{\phi\varphi}=W_\phi\,\ell_0\,H_\varphi+H_\phi\,\ell_0\,W_{\widetilde{\varphi}},
\end{align*}
where $\ell_0 : L^2(\mathbb{R_+})\rightarrow L^2(\mathbb{R})$ is the
zero extension operator. Additionally, from the last two identities,
it follows that
\begin{equation*}
W\!H_{\phi\varphi}=W_\phi\,\ell_0\,W\!H_\varphi+H_\phi\,\ell_0\,W\!H_{\widetilde{\varphi}}
\end{equation*}
and
\begin{equation}\label{eqseg}
W\!H_{\phi\varphi}=W\!H_\phi\,\ell_0\,W\!H_\varphi+H_\phi\,\ell_0\,W\!H_{\widetilde{\varphi}-\varphi}\;.
\end{equation}

Let $H^{\infty}(\mathbb{C_-})$ denote the set of all bounded and
analytic functions in $\mathbb{C_-}=\{z \in \mathbb{C}\,:\,
\mbox{Im}\,z < 0\}$, and let $H_{-}^{\infty}(\mathbb{R})$ be the set
of all functions in $L^{\infty}(\mathbb{R})$ that are non-tangential
limits of elements in $H^{\infty}(\mathbb{C_-}).$ We will also use
$H_{+}^{\infty}(\mathbb{R})$, which is defined in the obvious
corresponding way. Due to \eqref{eqseg}, if we consider $\phi \in
H_{-}^{\infty}(\mathbb{R})$ or $\varphi$ being an even function,
then we obtain a multiplicative relation
\begin{equation}\label{eq6}
W\!H_{\phi\varphi}=W\!H_\phi\,\ell_0\,W\!H_\varphi\;.
\end{equation}

Note that if the Fourier symbol of a Wiener-Hopf operator admits a
factorization of the form $\varphi_-\,\psi\,\varphi_+$, where
$\varphi_\pm \in H_{\pm}^{\infty}(\mathbb{R})$ and $\psi \in
L^{\infty}(\mathbb{R})$, it is possible to apply to the Wiener-Hopf
operator the multiplicative property
$$W_{\varphi_-\,\psi\,\varphi_+}=W_{\varphi_-}\,\ell_0\,W_\psi\,\ell_0\,W_{\varphi_+}$$
(see e.g. \cite [Proposition 2.17]{BoKaSp02}). With a convenient
change, it is possible to construct for Wiener-Hopf plus Hankel
operators a corresponding result as the one known for Wiener-Hopf
operators. In the present case we may apply the multiplicative
property on the left if the left factor belongs to
$H_{-}^{\infty}(\mathbb{R})$ and on the right if the right factor
is an even function, like the following proposition asserts.

\begin{proposition}\label{prop2}
Let $\varphi,\psi,\phi \in L^{\infty}(\mathbb{R})$. If $\varphi
\in H_{-}^{\infty}(\mathbb{R})$ and $\phi=\widetilde{\phi}$, then
the following operator factorization takes place:
\begin{align*}
W\!H_{\varphi\,\psi\,\phi}&=W\!H_\varphi\,\ell_0\,W\!H_\psi\,\ell_0\,W\!H_\phi\\
                        &=W_\varphi\,\ell_0\,W\!H_\psi\,\ell_0\,W\!H_\phi.
\end{align*}
\end{proposition}

\begin{proof}
From the hypothesis $\varphi \in H_{-}^{\infty}(\mathbb{R})$, we may
apply the already presented multiplicative relation for Wiener-Hopf
plus Hankel operators, see \eqref{eq6}. Thus
\begin{equation}\label{eq8}
W\!H_{\varphi\:\psi \phi}=W\!H_\varphi\,\ell_0\,W\!H_{\psi \phi}.
\end{equation}
In addition, since $\phi=\widetilde{\phi}$, it also follows from
\eqref{eq6} that
\begin{equation}\label{eq9}
W\!H_{\psi \phi}=W\!H_\psi\,\ell_0\,W\!H_\phi.
\end{equation}
From \eqref{eq8} and \eqref{eq9}, we have that
\begin{equation}\label{eq18}
W\!H_{\varphi\,\psi\,\phi}=W\!H_\varphi\,\ell_0\,W\!H_\psi\,\ell_0\,W\!H_\phi.
\end{equation}
Since $\varphi \in H_{-}^{\infty}(\mathbb{R})$, we have
$H_\varphi=0$ due to the structure of the Hankel operators.
Therefore $W\!H_\varphi=W_\varphi$ and it follows from \eqref{eq18}
that
$W\!H_{\varphi\,\psi\,\phi}=W_\varphi\,\ell_0\,W\!H_\psi\,\ell_0\,W\!H_\phi$.
\end{proof}

From (\ref{eq6}) we have that if the symbol of the Wiener-Hopf plus Hankel
operator is factorized in such a way that the right factor is an even
function, this leads to a factorization of the Wiener-Hopf plus Hankel
operator, where a Wiener Hopf plus Hankel operator with an even symbol
appears. Due to the multiplicative relation for Wiener-Hopf plus Hankel
operators (see (\ref{eq6})), we conclude that the Wiener Hopf plus Hankel
operator with an even symbol is an invertible operator. So, we end this
section with this result.

\begin{proposition}\label{prop1}
If $\phi_{e} \in \mathcal{G}L^{\infty}(\mathbb{R})$ and
$\widetilde{\phi_e}=\phi_e$, then $W\!H_{\phi_{e}}$ is invertible
and its inverse is the operator
$\ell_0\,W\!H_{\phi_{e}^{-1}}\,\ell_0 : L^2(\mathbb{R}_+)\rightarrow
L^2_+(\mathbb{R})$.
\end{proposition}

\begin{proof}
On the one hand, we have
\begin{equation}\label{eq7}
W\!H_{\phi_{e}\cdot
\phi_{e}^{-1}}\,\ell_0=W\!H_1\,\ell_0=W_1\,\ell_0=I_{L^2(\mathbb{R_+})},
\end{equation}
where $I_{L^2(\mathbb{R_+})}$ represents the identity operator in
$L^2(\mathbb{R_+})$. On the other hand, since $\phi_{e} \in
\mathcal{G}L^{\infty}(\mathbb{R})$ and $\widetilde{\phi_e}=\phi_e$,
then $\widetilde{\phi_e^{-1}}=\phi_e^{-1}$ and therefore we may
apply the \mbox{multiplicative} relation for Wiener-Hopf plus Hankel
operators. So we have
\begin{equation}\label{eq5}
W\!H_{\phi_{e}\cdot \phi_{e}^{-1}}
=W\!H_{\phi_{e}}\,\ell_0\,W\!H_{\phi_{e}^{-1}}.
\end{equation}
Thus, combining \eqref{eq7} and \eqref{eq5}, we get that
\begin{equation}\label{eqast1}
W\!H_{\phi_{e}}\,\ell_0\,W\!H_{\phi_{e}^{-1}}\,\ell_0=I_{L^2(\mathbb{R_+})}.
\end{equation}
In the same way, we obtain that
\begin{equation}\label{eqast2}
\ell_0\,W\!H_{\phi_{e}^{-1}}\,\ell_0\,W\!H_{\phi_{e}}=I_{L^2_+(\mathbb{R})}.
\end{equation}
Therefore, (\ref{eqast1})--(\ref{eqast2}) show that
$W\!H_{\phi_{e}}$ is invertible and its inverse is
$\ell_0\,W\!H_{\phi_{e}^{-1}}\ell_0$.
\end{proof}

\section{$AP$ Factorizations}\label{sect4}

We begin this section with the definition of a new kind of $AP$
factorization, the $AP$ asymmetric factorization. This definition was
initially motivated by the role of the so-called {\em $APW$ factorization}
in the theory of Wiener-Hopf operators with $APW$ Fourier symbols
\cite{BoKaSp02}, and by the recent works about Toeplitz plus Hankel
operators \cite{BaTo04, Eh04} and convolution type operators with
symmetry~\cite{CaSp05, CaSpTe04mn}. Additionally, it extends the
corresponding concept introduced in \cite{NoCa04} for the subclass of almost
periodic Fourier symbols $APW$.

\begin{definition}\label{def4.1}
We will say that a function $\phi \in \mathcal{G}AP$ admits an \textit{$AP$
asymmetric factorization} if it can be represented in the form
\begin{equation*}
\phi=\phi_{-} \; e_{\lambda} \; \phi_{e}
\end{equation*}
where $\lambda \in{\mathbb{R}} $, $e_{\lambda}(x)=e^{i\lambda x}$, $x \in
\mathbb{R}$, $\phi_{-} \in \mathcal{G}AP^{-}$, $\phi_{e} \in
\mathcal{G}L^\infty(\mathbb{R})$ with $\widetilde{\phi_e}=\phi_e$.

The particular case of an $AP$ asymmetric factorization with $\lambda=0$ will
be referred to as a \textit{canonical $AP$ asymmetric factorization}.
\end{definition}

\begin{figure}[h]
\hfil
\includegraphics[width=5.5cm]{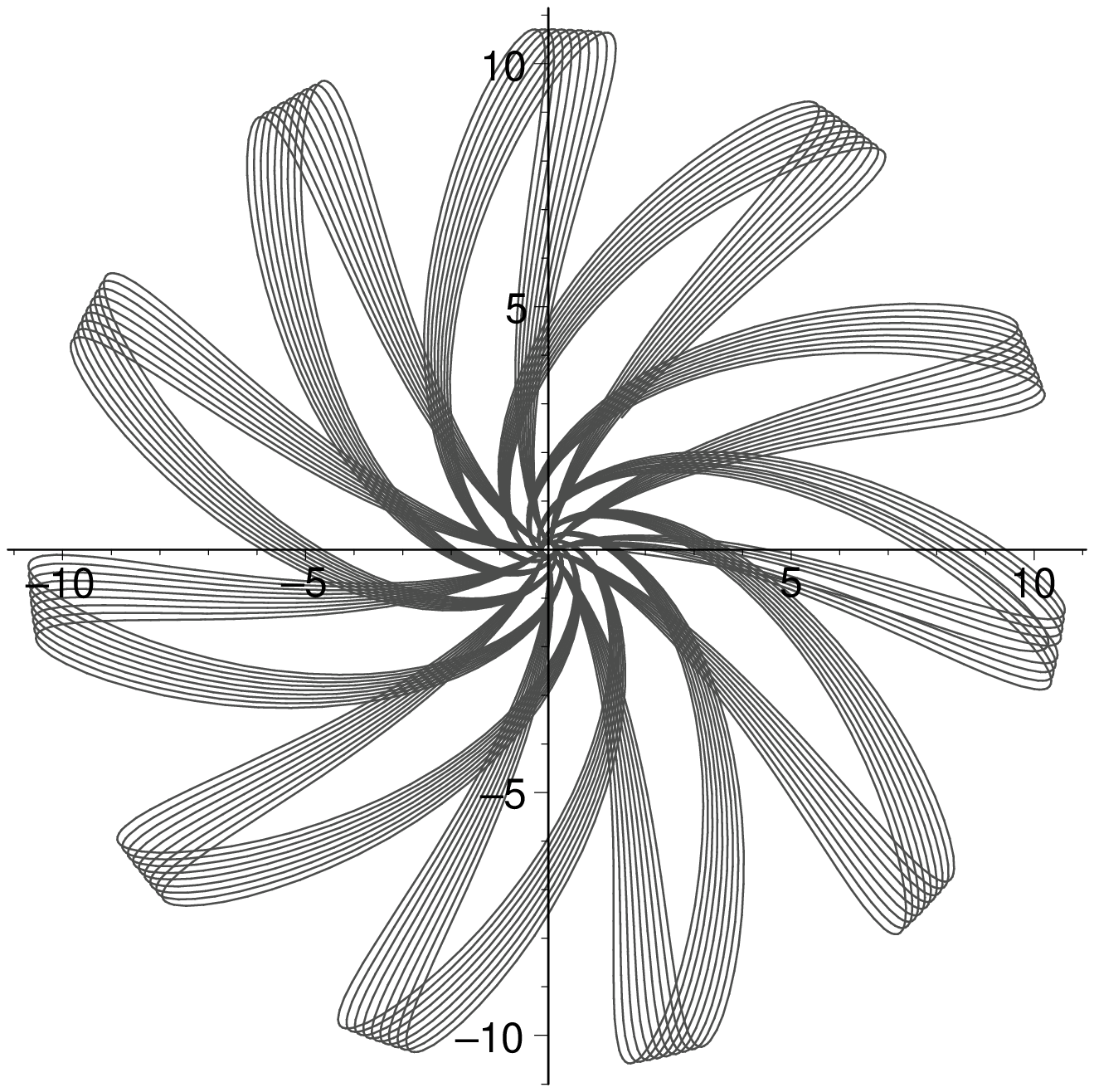}\hskip1cm\includegraphics[width=5.5cm]{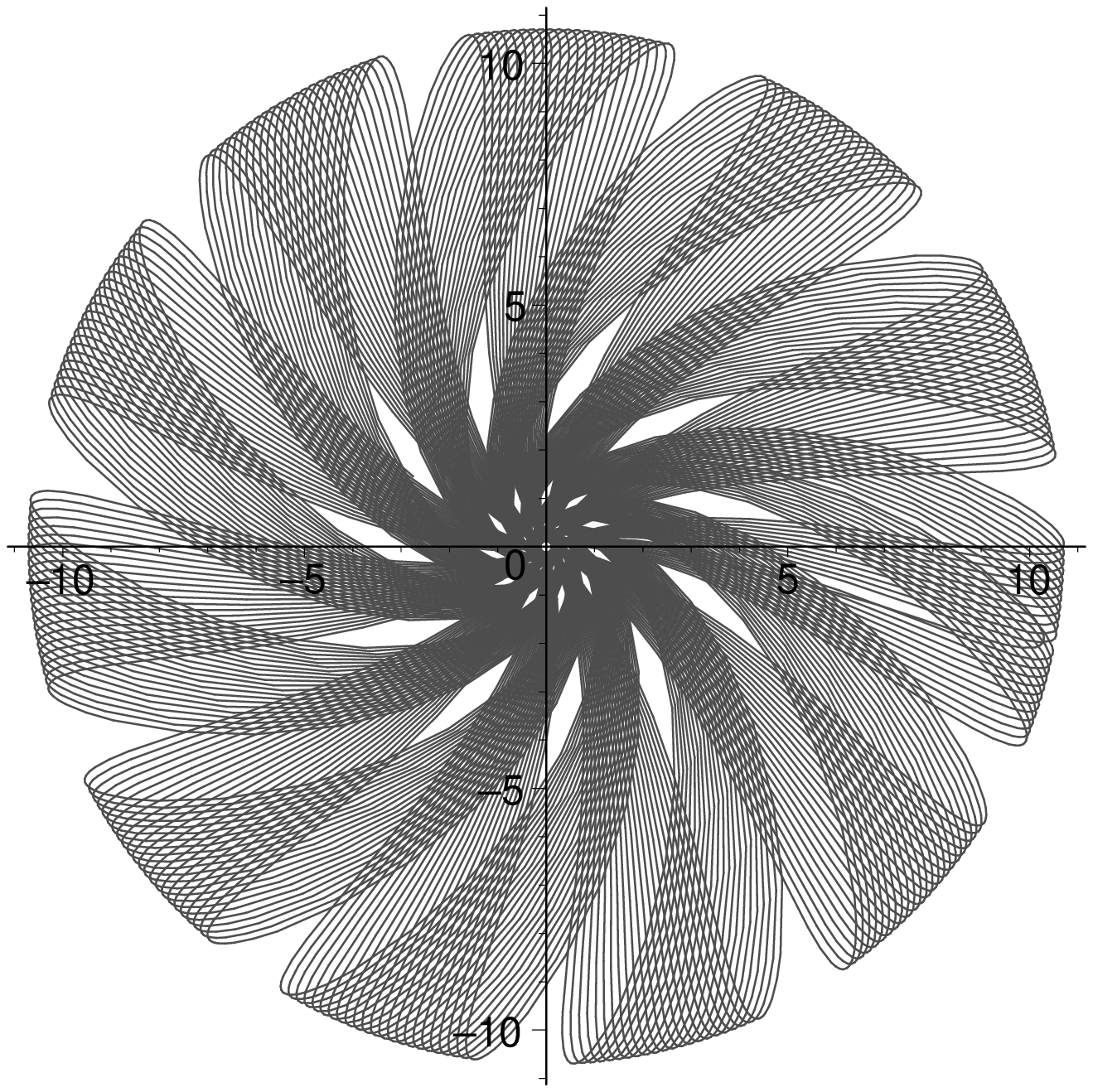}
\hfil
\includegraphics[width=12cm]{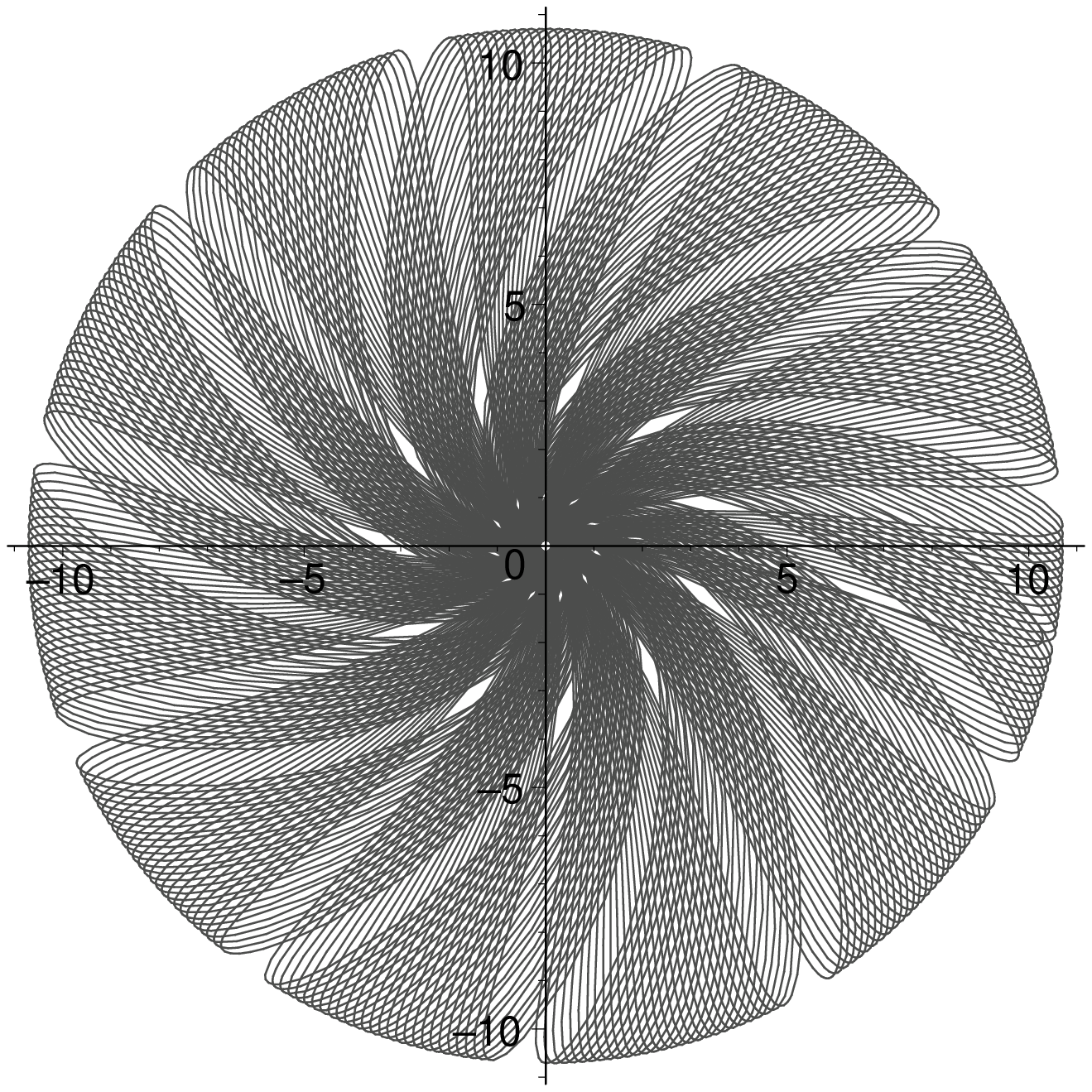}
\caption{The image of $\phi(x)=e^{-i2x+ e^{i\pi
x}}\ln((\arctan(100x^2)+\frac{\pi}{2})^{e^{-2\sin(\pi x)}})$ for
$x$ in $[-100,100]$, $[-250,250]$ and $[-300,300]$,
respectively.}\label{fig3}
\end{figure}

\begin{example}
Consider the function $\phi$ defined by
\begin{equation*}
\phi(x)=e^{\displaystyle{-i2x+ e^{\displaystyle{i\pi
x}}}}\,\ln\!\left(\bigg(\arctan(100x^2)+\frac{\pi}{2}\bigg)^{\displaystyle{e^{\displaystyle{-2\sin(\pi
x)}}}}\right),
\end{equation*}
for all $x \in \mathbb{R}$. Figure \ref{fig3} shows the image of
$\phi(x)$ (when $x$ belongs to three different intervals). We may
rewrite $\phi$ as
\begin{equation*}
\phi(x)=e^{\displaystyle{e^{-i\pi x}}}
e^{\displaystyle{-i2x}}\,\ln\!\bigg(\arctan(100x^2)+\frac{\pi}{2}\bigg),
\end{equation*}
for all $x \in \mathbb{R}$. Considering $\phi_-$ and $\phi_e$ given by
\begin{align*}
\phi_-(x)&=e^{\displaystyle{e^{-i\pi x}}},\\
\phi_e(x)&=\ln\bigg(\arctan(100x^2)+\frac{\pi}{2}\bigg),
\end{align*}
for all $x \in \mathbb{R}$, we have $\phi_- \in \mathcal{G}AP^-$ and
$\phi_{e} \in \mathcal{G}L^\infty(\mathbb{R})$ such that
$\widetilde{\phi_e}=\phi_e$. Therefore, it follows that $\phi$ admits an
$AP$ asymmetric factorization.

As for more examples, in $APW$ we find an endless number of functions which
have an $AP$ asymmetric factorization. Indeed, in \cite{NoCa04} it is proved
that every function in $APW$ admits an $APW$ asymmetric factorization. Since
$AP$ asymmetric factorization is a generalization of $APW$ asymmetric
factorization, it results that every function in $APW$ admits an $AP$
asymmetric factorization.
\end{example}

It is interesting to clarify that, when existing, the $AP$
asymmetric factorization of a function is unique up to a constant,
like it is stated in the following proposition.

\begin{proposition}
Let $\phi \in \mathcal{G}AP$. Suppose that $\phi$ admits two $AP$ asymmetric
\mbox{factorizations}:
\begin{align*}
\phi &= \phi_{-}^{(1)} \; e_{\lambda_1} \; \phi_{e}^{(1)}\;,\\
\phi &= \phi_{-}^{(2)} \; e_{\lambda_2} \; \phi_{e}^{(2)}\;.
\end{align*}
Then $\lambda_1=\lambda_2$, $\phi_{-}^{(1)}=\gamma \phi_{-}^{(2)}$ and
$\phi_{e}^{(1)}=\gamma^{-1} \phi_{e}^{(2)}$, $\gamma \in \mathbb{C}\setminus
\{0\}$.
\end{proposition}

\begin{proof}
The equality $\phi_{-}^{(1)} \; e_{\lambda_1} \; \phi_{e}^{(1)}=\phi_{-}^{(2)}
\; e_{\lambda_2} \; \phi_{e}^{(2)}$, implies that
\begin{equation}\label{eq11}
(\phi_{-}^{(2)})^{-1} \; \phi_{-}^{(1)} \; e_{\lambda_1} \; = e_{\lambda_2} \;
\phi_{e}^{(2)} \;(\phi_{e}^{(1)})^{-1}.
\end{equation}
Assume, without loss of generality, that $\lambda_1\leq\lambda_2$. Then
$\lambda=\lambda_1-\lambda_2\leq0$. From \eqref{eq11} it follows that
\begin{equation}\label{eq12}
(\phi_{-}^{(2)})^{-1} \; \phi_{-}^{(1)} \; e_{\lambda} \; =\phi_{e}^{(2)} \;
(\phi_{e}^{(1)})^{-1}.
\end{equation}
Since the right-hand side of \eqref{eq12} is an even function,
$(\phi_{-}^{(2)})^{-1} \; \phi_{-}^{(1)} \; e_{\lambda}$ is also an even
function. Put
\begin{equation}\label{eq13}
\varphi=(\phi_{-}^{(2)})^{-1} \; \phi_{-}^{(1)}.
\end{equation}
Thus
$\varphi(x)\; e_\lambda(x)=\widetilde{\varphi}(x)\;\widetilde{e_\lambda}(x)$,
i.e.
$\varphi(x)\;e_\lambda(x)=\widetilde{\varphi}(x)\;e_{-\lambda}(x)$,
or equivalently
\begin{eqnarray}
\varphi(x)\;e_{2\lambda}(x)=\widetilde{\varphi}(x).\label{eq14}
\end{eqnarray}
On the one hand, since $\varphi \in \mathcal{G}AP^{-}$, we may
apply the well-known characterization of $\mathcal{G}AP^-$ which
assures the existence of a $\psi \in AP^{-}$ such that
$\varphi=e^{\psi}$ (cf. e.g. \cite[Lemma 3.4]{BoKaSp02}). On the
other hand, because $\widetilde{\varphi} \in \mathcal{G}AP^{+}$,
by a corresponding characterization of $\mathcal{G}AP^+$, there
exists a $\eta \in AP^{+}$ such that
$\widetilde{\varphi}=e^{\eta}$. From \eqref{eq14}, it follows that
\begin{equation*}
e^{\psi(x)+i2\lambda x}=e^{\eta(x)},
\end{equation*}
which implies that $\lambda=0$ and $\psi \in AP^- \cap AP^+$, i.e.,
$\lambda_1=\lambda_2$ and $\psi$ is a complex constant function. From \eqref{eq13}, we
get $\phi_{-}^{(1)}=\gamma \phi_{-}^{(2)}$ with $\gamma \in
\mathbb{C}\backslash\{0\}$. By \eqref{eq12}, we obtain
$\phi_{e}^{(1)}=\gamma^{-1} \phi_{e}^{(2)}$.
\end{proof}

Let us recall that $\phi \in \mathcal{G}AP$ is said to admit a {\em
right $AP$ factorization}~\cite{BoKaSp02} if $\phi=\varphi_-
e_\lambda \varphi_+$, where $\varphi_{-} \in \mathcal{G}AP^-$,
$\varphi_{+} \in \mathcal{G}AP^+$, and $\lambda \in \mathbb{R}$. In
addition, if $\lambda =0$ this factorization is called a {\em
canonical right $AP$ factorization}. The $AP$ asymmetric
factorization is related to a special case of {\em right $AP$
factorization}, which we will call $AP$ antisymmetric factorization.
In this new kind of factorization a strong dependence between the
left and the right factor occurs (as we may realize in the next
definition).

\begin{definition}
A function $\phi \in \mathcal{G}AP$ admits an \textit{$AP$ antisymmetric
factorization} if it is possible to write
\begin{equation*}
\phi=\phi_{-} \; e_{2\lambda} \; \widetilde{\phi_{-}^{-1}}
\end{equation*}
where $\lambda \in{\mathbb{R}}$, $e_{2\lambda}(x)=e^{2 i\lambda x}$, $x \in
\mathbb{R}$, and $\phi_{-} \in \mathcal{G}AP^-$.
\end{definition}

The following proposition shows how the $AP$ asymmetric
factorization and the $AP$ antisymmetric factorization are related.
That is, $\phi$ has an $AP$ asymmetric factorization if and only if
$\Phi=\phi\widetilde{\phi^{-1}}$ has an $AP$ antisymmetric
factorization.

\begin{proposition}\label{prop3}
Let $\phi \in \mathcal{G}AP$ and put $\Phi=\phi\widetilde{\phi^{-1}}$.
\begin{itemize}
\item[{\rm(a)}] If $\phi$ admits an $AP$ asymmetric factorization,
$\phi=\phi_{-} \; e_{\lambda} \; \phi_{e}$, then $\Phi$ admits an $AP$
antisymmetric factorization with the same factor $\phi_-$ and the same index
$\lambda$. \item[{\rm(b)}] If $\Phi$ admits an $AP$ antisymmetric
factorization, $\Phi=\psi_{-} \; e_{2\lambda} \; \widetilde{\psi_{-}^{-1}}$,
then $\phi$ admits an $AP$ asymmetric factorization with the same minus
factor $\psi_-$, the same index $\lambda$ and the even factor
$\phi_e=e_{-\lambda}\psi_-^{-1}\phi$.
\end{itemize}
\end{proposition}

\begin{proof}
(a) From the $AP$ asymmetric factorization of $\phi$, $\phi=\phi_{-} \;
e_{\lambda} \; \phi_{e}$, we have
$$\widetilde{\phi^{-1}}=\phi_e^{-1} \;
e_{\lambda} \; \widetilde{\phi_-^{-1}}\,,$$ with $\phi_{-} \in
\mathcal{G}AP^-$. Hence
$$\Phi=\phi\widetilde{\phi^{-1}}=\phi_{-}
\; e_{2\lambda} \; \widetilde{\phi_{-}^{-1}}\,.$$

(b) It follows from the definition of the factor $\phi_e$ that
$\phi=\psi_{-} \; e_{\lambda} \; \phi_{e}$. Thus it remains to prove that
$\phi_e$ is an even function. Once again by the definition of $\phi_e$, we
obtain
\begin{equation*}
\widetilde{\phi_e}=e_{\lambda}\widetilde{\psi_-^{-1}}\widetilde{\phi}=
e_{\lambda}e_{-2\lambda}\psi_-^{-1}\phi=e_{-\lambda}\psi_-^{-1}\phi=\phi_e,
\end{equation*}
since $\widetilde{\psi_-^{-1}}\widetilde{\phi}
=e_{-2\lambda}\psi_-^{-1}\phi$ (due to the $AP$ antisymmetric factorization
of $\Phi$). Therefore $\phi_e$ is an even function.
\end{proof}

Since $AP$ antisymmetric factorization is a special case of right
$AP$ factorization and, by Proposition~\ref{prop3}, we already
know how to relate $AP$ antisymmetric factorization with $AP$
asymmetric factorization, it is natural to wonder how $AP$
asymmetric factorization and  right $AP$ factorization are
related. The next theorem gives the answer to this question.

\begin{theorem}\label{th2}
Let $\phi \in \mathcal{G}AP$. If $\phi$ admits a right $AP$ factorization,
$$\phi=\varphi_- e_\lambda \varphi_+\,,$$
then $\phi$ admits an $AP$ asymmetric factorization,
$$\phi=\phi_- e_\lambda
\phi_e\,,$$ with $\phi_-=\varphi_- \widetilde{\varphi_+^{-1}}$, and
$\phi_e=\widetilde{\varphi_+}\varphi_+$.
\end{theorem}

\begin{proof}Suppose that $\phi$ admits a right $AP$ factorization,
i.e., $\phi=\varphi_- e_\lambda \varphi_+$, where $\varphi_{-} \in
\mathcal{G}AP^-$, $\varphi_{+} \in \mathcal{G}AP^+$. Considering
$\Phi=\phi\widetilde{\phi^{-1}}$, we have
\begin{align}\label{eq20}
\Phi=\varphi_- e_\lambda \varphi_+ \widetilde{\varphi_+^{-1}} e_\lambda
\widetilde{\varphi_-^{-1}}
 =\varphi_- \widetilde{\varphi_+^{-1}} e_{2\lambda} \varphi_+ \widetilde{\varphi_-^{-1}}.
\end{align}
Since $\varphi_{-} \in \mathcal{G}AP^-$ and $\varphi_{+} \in
\mathcal{G}AP^+$, then $\widetilde{\varphi_+^{-1}}\in \mathcal{G}AP^-$,
$\widetilde{\varphi_-^{-1}} \in \mathcal{G}AP^+$ and therefore $\varphi_-
\widetilde{\varphi_+^{-1}}\in \mathcal{G}AP^-$ and $\varphi_+
\widetilde{\varphi_-^{-1}}\in \mathcal{G}AP^+$. Putting $\phi_-=\varphi_-
\widetilde{\varphi_+^{-1}}$, it follows from \eqref{eq20} that
$$\Phi=\phi_- e_{2\lambda} \widetilde{\phi_-^{-1}}\,.$$
Since $\phi_{-} \in \mathcal{G}AP^-$, it results that $\Phi$ admits a $AP$
antisymmetric factorization. By Proposition \ref{prop3}, that implies that
$\phi$ admits a $AP$ asymmetric factorization, $\phi=\phi_- e_\lambda
\phi_e$, with $\phi_e= e_{-\lambda} \phi_-^{-1}\phi$. Rewriting $\phi_-$ and
$\phi_e$ by using the factors of the right $AP$ factorization, $\varphi_-$
and $\varphi_+$, we have
$$\phi_-=\varphi_-
\widetilde{\varphi_+^{-1}}\,,\qquad
\phi_e=\widetilde{\varphi_+}\varphi_+\,.$$
\end{proof}

\begin{corollary}\label{cor1}
Let $\phi \in \mathcal{G}AP$. If $\phi$ admits a canonical right $AP$
factorization,
$$\phi=\varphi_- \varphi_+\,,$$
then $\phi$ admits a canonical $AP$ asymmetric factorization,
$$\phi=\phi_-
\phi_e\,,$$ with $\phi_-=\varphi_- \widetilde{\varphi_+^{-1}}$, and
$\phi_e=\widetilde{\varphi_+}\varphi_+$.
\end{corollary}

\begin{proof}
The result is a direct consequence of Theorem~\ref{th2}, if we
take $\lambda = 0$.
\end{proof}

\section{Main Results: Invertibility of Wiener-Hopf plus Hankel
Operators}\label{sect5}

In this section, we present an invertibility and Fredholm
criterion, as well as an explicit formula for the (one-sided and
two-sided) inverses of Wiener-Hopf plus Hankel operators with
almost periodic Fourier symbols. Both results are obtained in
terms of an $AP$ asymmetric factorization of the Fourier symbol of
the operators in study.

\begin{theorem}\label{th1}
Let $\phi \in \mathcal{G}AP$ admit an $AP$ asymmetric factorization
$\phi=\phi_{-} \; e_{\lambda} \; \phi_{e}$.
\begin{itemize}
\item[{\rm(a)}] If $\lambda < 0$, then $W\!H_\phi$ is properly
right-Fredholm and right-invertible. \item[{\rm(b)}] If $\lambda > 0$, then
$W\!H_\phi$ is properly left-Fredholm and left-invertible.
\item[{\rm(c)}] If $\lambda = 0$, then $W\!H_\phi$ is invertible.
\end{itemize}
\end{theorem}

\begin{proof}
In the case where $\lambda < 0$, we have that $e_{\lambda} \in AP^{-}$. Since
$AP^{-}=AP \, \cap \, H_{-}^{\infty}(\mathbb{R})$, it holds that
$e_{\lambda}\in H_{-}^{\infty}(\mathbb{R})$ and hence
\begin{equation}\label{eqrel1}
W\!H_\phi=W_{\phi_{-}}\,\ell_0\,W_{e_{\lambda}}\,\ell_0\,W\!H_{\phi_{e}},
\end{equation}
due to Proposition~\ref{prop2} and also taking into account that, because
$e_{\lambda}\in H_{-}^{\infty}(\mathbb{R})$,
$W\!H_{e_{\lambda}}=W_{e_{\lambda}}$. Since $\phi_{-} \in
\mathcal{G}AP^{-}$, by the characterization of $\mathcal{G}AP^-$, there
exists a $\psi \in AP^{-}$ such that $\phi_{-}=e^{\psi}$. Thus, the mean
motion of $\phi_{-}$ is zero and due to the {\em
Gohberg-Feldman/Coburn-Douglas Theorem} (stated in the first section),
$W_{\phi_{-}}$ is invertible. From Proposition \ref{prop1}, we know that
$W\!H_{\phi_{e}}$ is invertible. Therefore, since
$\ell_0:L^2(\mathbb{R_+})\rightarrow L^2_+(\mathbb{R})$ is also an
invertible operator, (\ref{eqrel1}) shows that $W\!H_\phi$ is equivalent to
$W_{e_{\lambda}}$. Once again, by the Theorem of {\em
Gohberg-Feldman/Coburn-Douglas}, since the mean motion of ${e_{\lambda}}$ is
$\lambda < 0$, we have that the operator $W_{e_{\lambda}}$ is properly right
Fredholm and right-invertible. Consequently, due to the equivalence relation
(\ref{eqrel1}), the operator $W\!H_\phi$ is also properly right-Fredholm and
right-invertible. This completes the proof of part (a).

Part (b) can be derived from part (a) by passage to adjoint operators.

Finally, let us now suppose that $\lambda = 0$. Then
$\phi=\phi_{-}\: \phi_{e}$ and
$W\!H_\phi=W_{\phi_{-}}\,\ell_0\,W\!H_{\phi_{e}}$. Since
$W_{\phi_{-}}$ and $W\!H_{\phi_{e}}$ are invertible, then
$W\!H_\phi$ is also invertible.
\end{proof}

In order to proceed to the next result, we recall here the concept of
\emph{reflexive generalized invertibility}. Let $T:X\rightarrow Y$ be a
bounded linear operator acting between Banach spaces. $T$ is said to be
\emph{reflexive generalized invertible} if there exists a bounded linear
operator $T^-:Y\rightarrow X$ such that $TT^-T=T$ and $T^-TT^-=T^-$. In this
case, the operator $T^-$ is referred to as the \emph{reflexive generalized
inverse} of $T$. A linear bounded one-sided or two-sided invertible operator
is also a reflexive generalized invertible operator.

After having reached an invertibility criterion for the
Wiener-Hopf plus Hankel operators, we start looking for an
explicit formula for the reflexive generalized inverses of these
Wiener-Hopf plus Hankel operators. The obtained formula is
expressed in terms of the factors of the $AP$ asymmetric
factorization. From the value of $\lambda$ of the middle factor of
the $AP$ asymmetric factorization, $e_\lambda$, it is possible to
distinguish that the reflexive generalized inverse is in fact a
right-inverse, left-inverse or inverse. I.e., depending on
$\lambda < 0$, $\lambda
> 0$ and $\lambda = 0$, we obtain the right-inverse, the
left-inverse and the inverse of ${W\!H}_\phi$, respectively. As we
will see in the theorem below, the explicit formula for the
reflexive generalized inverse of  ${W\!H}_\phi$ is given in terms of
an arbitrary extension operator $\ell:L^2(\mathbb{R_+})\rightarrow
L^2(\mathbb{R})$. This means that the reflexive generalized inverse
of ${W\!H}_\phi$ is independent of the choice of $\ell$ (and
therefore several choices are allowed, like for instance $\ell =
\ell_0$ or $\ell = \ell^e$).

\begin{theorem}\label{th3}
If $\phi \in \mathcal{G}AP$ admits an $AP$ asymmetric factorization
$$\phi=\phi_{-} \; e_{\lambda} \; \phi_{e}\,,$$ then we obtain a reflexive
generalized inverse of $W\!H_\phi$ defined by
\begin{equation*}
{W\!H}_\phi^- = \ell_0 r_+ \mathcal F^{-1}\phi^{-1}_{e}\cdot \mathcal F \ell^e
r_+ {\mathcal F}^{-1}e_{-\lambda}\cdot \mathcal F \ell^e r_+ \mathcal
F^{-1}\phi^{-1}_{-}\cdot \mathcal F
 \ell:L^2(\mathbb{R_+})\rightarrow L^2_+(\mathbb{R}),
\end{equation*}
where $\ell:L^2(\mathbb{R_+})\rightarrow L^2(\mathbb{R})$ denotes an arbitrary
extension operator.

Additionally, in a more detailed way:
\begin{itemize}
\item[{\rm(a)}] if $\lambda < 0$, then ${W\!H}_\phi^-$ is the right-inverse of
$W\!H_\phi$; \item[{\rm(b)}] if $\lambda > 0$, then ${W\!H}_\phi^-$ is the
left-inverse of $W\!H_\phi$; \item[{\rm(c)}] if $\lambda = 0$, then
${W\!H}_\phi^-$ is the inverse of $W\!H_\phi$.
\end{itemize}
\end{theorem}

\begin{proof}
From the $AP$ asymmetric factorization $\phi=\phi_{-} \: e_{\lambda} \:
\phi_{e}$, it directly follows that
\begin{equation*}
W\!H_\phi= r_{+}\,A_{-}\,E\,A_{e}\,\ell^{e}r_{+},
\end{equation*}
where $A_{-}=\mathcal F^{-1}\phi_{-}\cdot \mathcal F$, $E=\mathcal
F^{-1}e_{\lambda}\cdot \mathcal F$ and $A_{e}=\mathcal F^{-1}\phi_{e}\cdot
\mathcal F$.

(i) If $\lambda\leq0$, consider
\begin{align}\label{eq15x}
W\!H_\phi\, {W\!H}_\phi^-  &= r_{+}\,A_{-}\,E\,A_{e}\,\ell^{e} r_{+}\,\;
\ell_0 r_+ A_{e}^{-1} \ell^{e} r_{+} \, E^{-1} \ell^{e} r_{+} \;
A_{-}^{-1} \ell \notag\\
&= r_{+}\,A_{-}\,E\,A_{e}\,\ell^{e} r_{+} A_{e}^{-1} \ell^{e}
r_{+} \, E^{-1} \ell^{e} r_{+} \, A_{-}^{-1} \ell,
\end{align}
where the term $\ell_0 r_+$ was omitted due to the fact that $r_{+}\ell_0 r_+
=r_+$. Since $A_e^{-1}$ preserves the even property of its symbol, we may also
drop the first  $\ell^{e} r_{+}$ term in \eqref{eq15x}, and obtain
\begin{equation}\label{eq16a}
W\!H_\phi\, {W\!H}_\phi^- = r_{+}\,A_{-}\,E\, \ell^{e} r_{+} \, E^{-1}
\ell^{e} r_{+} \, A_{-}^{-1} \ell.
\end{equation}
Additionally, since in the present case (due to $\lambda\leq0$) $E^{-1}$ is a
{\em plus type factor} \cite{CaSpTe04mn, Sp85},  we have $\ell^{e} r_{+} \,
E^{-1} \ell^{e} r_{+}= E^{-1} \ell^{e} r_{+}$; also because $A_{-}$ is a {\em
minus type factor} it follows
\begin{equation}\label{eq16b}
W\!H_\phi\, {W\!H}_\phi^- = r_{+}\,A_{-}\ell^{e} r_{+} \, A_{-}^{-1} \ell =
r_{+}\, \ell = I_{L^2(\mathbb{R_+})},
\end{equation}
and we can directly realize that such identities do not depend on the
particular choice of the extension operator $\ell$.

(ii) If $\lambda\geq0$, we will now analyze the composition
\begin{align}\label{eq15c}
{W\!H}_\phi^- \,W\!H_\phi&=  \ell_0 r_+ A_{e}^{-1} \ell^{e} r_{+} \, E^{-1}
\ell^{e} r_{+} \; A_{-}^{-1} \ell \,\;
 r_{+}\,A_{-}\,E\,A_{e}\,\ell^{e} r_{+}\,.
\end{align}
In the present case $E^{-1}$ is a {\em minus type factor} and for this reason
 $\ell^{e} r_{+} E^{-1} \ell^{e} r_{+}$ $= \ell^{e} r_{+} \, E^{-1}$. The
same reasoning applies to the factor $A_{-}^{-1}$, and therefore the equality
(\ref{eq15c}) takes the form
\begin{align}\label{eq15f}
{W\!H}_\phi^- \,W\!H_\phi&=  \ell_0 r_+ A_{e}^{-1} \ell^{e} r_{+} \,
\,A_{e}\,\ell^{e} r_{+} = \ell_0 r_+ \,\ell^{e} r_{+} = \ell_0 r_+ =
I_{L^2_+(\mathbb{R})}\,,
\end{align}
where we have used the identity $\ell^e r_+ A_{e} \ell^e r_+ = A_{e} \ell^e
r_+$.

(iii) Intersecting the last two cases, (i) and (ii), it follows that for
$\lambda = 0$, the operator ${W\!H}_\phi^-$ is the (both-sided) inverse of
$W\!H_\phi$ (cf. (\ref{eq16b}) and (\ref{eq15f})).
\end{proof}

As a consequence of Theorem~\ref{th1} and of the relation between
right $AP$ factorizations and $AP$ asymmetric factorizations
presented in Theorem~\ref{th2} and Corollary~\ref{cor1}, we end up
with a curious result on the dependence between the invertibility of
Wiener-Hopf and Wiener-Hopf plus Hankel operators with the same $AP$
Fourier symbol. That is, from the invertibility of the  Wiener-Hopf
operator, we obtain the invertibility of the Wiener-Hopf plus Hankel
operator.

\begin{corollary}\label{Cor5.3}
Let $\phi \in \mathcal{G}AP$. If $W_\phi$ is invertible with $\phi$ having a
canonical right AP factorization,  then $W\!H_\phi$ is invertible.
\end{corollary}

\begin{proof}
Suppose that $\phi=\varphi_- \varphi_+$ is a canonical right $AP$ factorization
of $\phi$. By Corollary~\ref{cor1}, $\phi$ admits a canonical $AP$ asymmetric
factorization, $\phi=\phi_-\phi_e$, where
$$\phi_- = \varphi_- \,\,\widetilde{\varphi_+^{-1}}\,,\qquad \phi_e
= \widetilde{\varphi_+}\,\,\varphi_+\,.$$
From Theorem \ref{th1},
it follows that $W\!H_\phi$ is an invertible operator.
\end{proof}

\noindent{\bf Acknowledgements}: The work was supported by {\em
Funda\c{c}\~{a}o para a Ci\^{e}ncia e a Tecnologia} through {\em Unidade de
Investiga\c{c}\~{a}o Matem\'{a}tica e Aplica\c{c}\~{o}es} of University of Aveiro,
Portugal.

\end{document}